\documentclass[leqno,11pt]{amsart}

\usepackage{amssymb}
\usepackage{amsmath}
\usepackage{enumerate}
\usepackage[dvips]{graphicx}

\def\r{\mathbb{R}}
\def\R{\mathbb{R}}
\def\z{\mathbb{Z}}

\def\c{\mathbb{C}}
\def\C{\mathbb{C}}

\newtheorem{definition}{Definition}
\newtheorem{remark}{Remark}
\newtheorem{theorem}{Theorem}
\newtheorem{proposition}{Proposition}
\newtheorem{corollary}{Corollary}
\newtheorem{lemma}{Lemma}

\DeclareMathOperator{\arctanh}{arctanh}

\begin{document}

\title[Lagrangian translating solitons for MCF]{Translating solitons for Lagrangian mean curvature flow in complex Euclidean plane}

\author{Ildefonso Castro}
\address{Departamento de Matem\'{a}ticas \\
Universidad de Ja\'{e}n \\
23071 Ja\'{e}n, SPAIN} \email{icastro@ujaen.es}

\author{Ana M.~Lerma}
\address{Departamento de Matem\'{a}ticas \\
Universidad de Ja\'{e}n \\
23071 Ja\'{e}n, SPAIN} \email{alerma@ujaen.es}

\thanks{Research partially supported by a MEC-Feder grant MTM2007-61775}



\date{}

\begin{abstract}
Using certain solutions of the curve shortening flow, including
self-shrinking and self-expanding curves or spirals, we construct
and characterize many new examples of translating solitons for
mean curvature flow in complex Euclidean plane. They generalize
the Joyce, Lee and Tsui ones \cite{JLT} in dimension two. The
simplest (non trivial) example in our family is characterized as
the only (non totally geodesic) Hamiltonian stationary Lagrangian
translating soliton for mean curvature flow in complex Euclidean
plane.
\end{abstract}

\maketitle

\section{Introduction}

The mean curvature flow (in short MCF) of an immersion $\phi:
M\rightarrow \r^4$ of a smooth surface $M$ is a family of
immersions $F:M \times [0,\epsilon) \rightarrow \r^4$ parametrized
by $t$ that satisfies
\begin{equation}\label{MCF}
\frac{d}{dt}F_t(p)=H(p,t), \quad F_0=\phi,
\end{equation}
where $H(p,t)$ is the mean curvature vector of $F_t(M)$ at
$F_t(p)=F(p,t)$. The evolution of a Lagrangian surface in complex
Euclidean plane $\c^2$ by its mean curvature preserves its
Lagrangian character and it is called the {\em Lagrangian mean
curvature flow}.

Some interesting problems rather far from trivial in this setting
are, on the one hand, to understand the possible singularities
that can occur during the flow in finite time and, on the other
hand, if it is possible to show that the singularities for
Lagrangian MCF are isolated. A.\ Neves constructed in \cite{N}
examples of Lagrangians in $\c^2$ having the Lagrangian angle as
small as desired and for which the Lagrangian MCF develops a
finite-time singularity. But he also proved in \cite{N} that
assuming certain properties on the initial Lagrangian surface,
like almost calibrated, i.e.\ the oscillation of the Lagrangian
angle to be strictly smaller than $\pi$, if one rescales the flow
around a fixed point in space-time, connected components of this
rescaled flow converge to an area-minimizing union of planes.

In geometric flows such as the Ricci flow or the Lagrangian MCF,
singularities are often locally modelled on soliton solutions,
such as Lagrangians which are moved by rescaling or translation by
MCF. When the evolution is a homotethy we get the {\em
self-similar solutions} for MCF. In the Lagrangian context they
have been considered by several authors; see for example
\cite{CL}, \cite{JLT}, \cite{LW1} and \cite{LW2}. The study of
this type of solutions is hoped to give a better understanding of
the flow at a singularity since by Huisken's monotonicity formula
\cite{Hu}, any central blow-up of a finite-time singularity of the
mean curvature flow is a self-similar solution.

J.\ Chen and J.\ Li \cite{ChL} and M.-T.\ Wang \cite{Wa} proved
independently that there is no Type I singularity along the almost
calibrated Lagrangian mean curvature flow. Therefore it is of
great interest to understand dilations of the flow where the point
at which we center the dilation changes with the scale, called
Type II dilations, which converge to an eternal solution with
second fundamental form uniformly bounded. One of the most
important examples of Type II singularities is a class of eternal
solutions known as {\em translating solitons}, which are surfaces
which evolve by translating in space with constant velocity.

The eternal solution $F_t$, $t\in \r$, defined by
\begin{equation}\label{grim}
F_t(x,y)=(-\log\cos y + t,y,x,0), \
-\frac{\pi}{2}<y<\frac{\pi}{2}, \, x\in\r
\end{equation}
is called the {\em grim-reaper} and it is probably the most known
example of translating solution to MCF.

In \cite{NT}, A.\ Neves and G.\ Tian gave conditions that exclude
the existence of nontrivial translating solutions to Lagrangian
MCF. More precisely, they proved that translating solutions with
an $L^2$ bound on the mean curvature vector are planes and almost
calibrated translating solutions which are static are also planes.

D.~Joyce, Y.-I.~Lee and M.-P.~Tsui found out in \cite{JLT} new
surprising translating solitons for Lagrangian MCF with
oscillation of the Lagrangian angle arbitrarily small. They play
the same role as cigar solitons in Ricci flow and are important in
studying the regularity of Lagrangian MCF. Moreover, joint to the
grim-reaper (\ref{grim}), these examples show that the geometric
conditions on the above results in \cite{NT} are optimal.

In Section 2 we describe the main geometric properties of the
Lagrangian translating solitons and recall some examples.  Some
other interesting properties of them are studied in \cite{HanLi}.

In Section 3 we generalize Joyce, Lee and Tsui examples to a
considerable extent: It is remarked in \cite{NT} that they are
associated to planar curves $w$ in $\c$ such that
$w_t:=\sqrt{2t}w$, for $t>0$, is a solution to curve shortening
flow in $\c$. However, our general construction is based in two
families of planar curves $\alpha $ and $\omega $ depending on an
angular parameter $\varphi \in [0,\pi)$ (see Proposition 2) that
are special solutions to curve shortening flow (see Lemma 1), in
the sense that their flows are a kind of composition of dilations
and rotations with suitable velocities depending on $\varphi$. For
instance, in the case $\varphi =\pi / 2$ we must consider $\alpha
$ and $\omega $ spirals (i.e. travelling waves in the polar angle,
see \cite{CZ}) with opposite velocities; and in the case $\varphi
=0$, we require this time self-similar solutions for the curve
shortening flow of opposite characters, that is, $\alpha $ must be
a self-shrinking curve while $\omega$ must be a self-expanding
one. Just when in this particular case $\varphi =0$ we consider
$\alpha $ as a straight line passing through the origin, we arrive
at the above Joyce, Lee and Tsui examples (see Corollary 1).

In \cite{CL}, the authors classified the Hamiltonian stationary
Lagrangian self-similar solutions for Lagrangian mean curvature
flow in complex Euclidean plane. Three one-parameter families of
surfaces with different topologies (including embedded nontrivial
planes) appeared. In Section 4
 we characterize locally all our examples (see Theorem 1) in terms of
 an analytical condition on the Hermitian product of the position vector
 of the immersion and the translating vector that allow us
 {\em separation of variables}. As a
 consequence we get in Corollary 3 the classification of the Hamiltonian
 stationary Lagrangian translating solitons for Lagrangian mean curvature
flow in complex Euclidean plane. In contrast to the self-similar
case, only one example appears (see Corollary 2): a embedded
complete nontrivial plane given by
$$ \mathcal{M} = \{ (z,w)\in \c^2 \, : \, w^2=2 {\rm Re}z \, e^{-2i{\rm Im}z}, \,\, {\rm Re}z \geq 0 \}. $$
It corresponds in our construction to the simplest nontrivial
possible election of $\alpha $  (the circle $\alpha (t)=e^{it}$)
and $\omega$ (the line $\omega (s)=s$) in the particular case
$\varphi = 0$.

Joyce, Lee and Tsui examples are the only ones in our family with
oscillation of the Lagrangian angle arbitrarily small. Therefore
it would be very important to solve the open question if they can
arise as blow-ups of finite time singularities for Lagrangian mean
curvature flow.



\section{Preliminaries}

\subsection{Lagrangian surfaces in complex Euclidean plane}

In the complex Euclidean plane $\c^2$ we consider the bilinear Hermitian product defined by
\[
(z,w)=z_1\bar{w}_1+z_2\bar{w}_2, \quad  z,w\in\c^2.
\]
Then $\langle\, \, , \, \rangle = {\rm Re} (\,\, , \,)$ is the Euclidean metric on $\c^2$ and
$\omega = -{\rm Im} (\,,)$ is the Kaehler two-form given by $\omega (\,\cdot\, ,\,\cdot\,)=\langle
J\cdot,\cdot\rangle$, where $J$ is the complex structure on $\c^2$. We also consider the closed
complex-valued 2-form given by $\Omega = dz_1 \wedge dz_2$ and the Liouville 1-form $\lambda$ given
by $d\lambda = 2 \omega$.

Let $\phi:M \rightarrow \c^2$ be an isometric immersion of a surface $M$ into $\c^2$. $\phi $ is
said to be Lagrangian if  $\phi^* \omega = 0$.  Then we have $\phi^* T\c^2 =\phi_* TM \oplus J
\phi_* T M$, where $TM$ is the tangent bundle of $M$. The second fundamental form $\sigma $ of
$\phi $ is given by $\sigma (v,w)=JA_{Jv}w$, where $A$ is the shape operator, and so the trilinear
form $C(\cdot,\cdot,\cdot)=\langle \sigma(\cdot,\cdot), J \cdot \rangle $ is fully symmetric.

If $M$ is orientable and $\omega_M$ denotes the area form of $M$,
then $\phi^* \Omega = e^{i\beta} \omega_M$, where
$\beta:M\rightarrow \r /2\pi \z$ is called the {\em Lagrangian
angle} map of $\phi$ (see \cite{HL}). In general $\beta $ is a
multivalued function; nevertheless $d\beta $ is a well defined
closed 1-form on $M$ and its cohomology class is called the {\em
Maslov class}. When $\beta $ is a single valued function the
Lagrangian is called {\em zero-Maslov class} and if $\cos\beta
\geq \epsilon$ for some $\epsilon >0$ then the Lagrangian is said
to be {\em almost calibrated}.
 It is remarkable that $\beta$
satisfies (see for example \cite{SW})
\begin{equation}\label{beta}
J\nabla\beta=H=\Delta \phi,
\end{equation}
where $H$ is the mean curvature vector of $\phi$, defined by $H=
{\rm trace} \, \sigma$ and $\Delta $ is the Laplace operator of
the induced metric on $M$.

If $\beta $ is constant, say $\beta\equiv \beta_0$ or, equivalently $H=0$, then the Lagrangian
immersion $\phi $ is calibrated by ${\rm Re} (e^{-i\beta_0}\Omega)$ and hence area-minimizing. They
are referred as being {\em Special Lagrangian}.

A Lagrangian submanifold is called {\em Hamiltonian stationary} if
the Lagrangian angle $\beta $ is harmonic, i.e. $\Delta \beta =0$,
where $\Delta $ is the Laplace operator on $M$. Hamiltonian
stationary Lagrangian (in short HSL) surfaces are critical points
of the area functional with respect to a special class of
infinitesimal variations preserving the Lagrangian constraint;
namely, the class of compactly supported Hamiltonian vector fields
(see \cite{O}). Examples of HSL surfaces in $\c^2$ can be found in
\cite{A1}, \cite{CU2} and\cite{HR1}.

\subsection{Translating solitons for the mean curvature flow}

Let $\phi: M\rightarrow \r^4$ be an immersion of a smooth surface
$M$ in Euclidean 4-space. In geometric flows such as the Ricci
flow or the MCF, singularities are often locally modelled on
soliton solutions. In the case of MCF, one type of soliton
solutions of great interest are those moved by translating in the
Euclidean space. We recall that they must be of the following
form:

\begin{definition}\label{def}
An immersion $\phi: M\rightarrow \r^4$ is called a translating
soliton for mean curvature flow if
\begin{equation}\label{trl}
H={\bf e}^\perp
\end{equation}
for some nonzero constant vector ${\bf e}\in \R^4$, where ${\bf
e}^\perp $ denotes the normal projection of the vector $\bf e $
and $H$ is the mean curvature vector of $\phi $. The 1-parameter
family $F_t:=\phi + t {\bf e} $, $t\in\r$, is then solution of
(\ref{MCF}) and we call $\bf e$ a translating vector.
\end{definition}

 Any
translating soliton for MCF must be a gradient soliton, that is,
${\bf e}^\top=\nabla f $, for some smooth function $f:M\rightarrow
\r$, where ${\bf e}^\top $ denotes the tangent projection of the
vector $\bf e $. In fact, it is proved in \cite{JLT} that ${\bf
e}^\top=\nabla \langle \phi, {\bf e}\rangle $.

For Lagrangian translating solitons for MCF we point out the
following properties.

\begin{proposition}\label{properties}
Let $\phi: M \rightarrow \c^2$ be a Lagrangian translating soliton
for mean curvature flow with translating vector $\bf e$ and
Lagrangian angle map $\beta$. Then:
\begin{enumerate}
\item $ \beta = - \langle \phi , J{\bf e} \rangle + \beta_0 $, for
some constant $\beta_0$;
\item $\Delta \beta + \langle \nabla \beta , {\bf e} \rangle =0$;
\item $\Delta \langle \phi , {\bf e} \rangle= |H|^2$.
\end{enumerate}
\end{proposition}

\begin{proof}
Using (\ref{beta}) and (\ref{trl}) we have that $\nabla \beta =-(J {\bf e} )^\top$ and so $\nabla
\beta + \nabla \langle \phi , J{\bf e} \rangle =0$, which proves part 1. In addition, using
(\ref{beta}) again, $\Delta \beta =- \Delta  \langle \phi , J{\bf e} \rangle = - \langle \nabla
\beta,  {\bf e} \rangle$, which gives part 2. Finally, from (\ref{beta}) and (\ref{trl}) we deduce
$\Delta \langle \phi , {\bf e} \rangle = \langle H, {\bf e}^\perp \rangle = |H|^2$, which is part
3.
\end{proof}

In particular, part 1 in Proposition 1 says that a Lagrangian
translating soliton for MCF is always zero-Maslov class and from
part 3 we easily deduce that there are no compact Lagrangian
translating solitons for MCF.

By scaling and choosing a suitable coordinate system in
$\r^4\equiv \c^2$, we can assume that ${\bf e}=(1,0,0,0)\equiv
(1,0)\in \C^2$ without loss of generality.

\subsection{Examples of Lagrangian translating solitons}

The simplest examples of Lagrangian surfaces in $\C^2$ are usually
found as product of planar curves. If we look for translating
solitons for MCF in this family, we note that the {\em
grim-reaper} $F_t$, $t\in \r$, defined in (\ref{grim}) can be
written as
\[
F_t(x,y)= (\gamma (y),x) + t (1,0), \ \gamma(y)=-\log\cos y + i \,
y, \ -\frac{\pi}{2}<y<\frac{\pi}{2}, \, x\in\r,
\]
so $\gamma $ being the graph of $-\log\cos y$ that we will call
the grim-reaper curve. We can parameterize $\gamma $ by arc length
$s= 2 \arctanh (\tan y/2) $ obtaining
\begin{equation}\label{grcurve}
\gamma (s)=
(\log \cosh s, 2 \arctan (\tanh s /2)), \ s\in \R .
\end{equation}
It is remarkable that the
curvature $\kappa_\gamma$ of $\gamma $ verifies $\kappa_\gamma
(s)=-\gamma_2'(s)=-1/\cosh s=-1/e^{\gamma_1(s)}$.

Using precisely this last property, it is an exercise to check that the product immersion
\begin{equation}\label{grimprod}
(s_1,s_2)\in \R^2 \longrightarrow (\gamma(s_1),\gamma(s_2))\in
\c^2
\end{equation}
is a translating soliton for MCF with translating vector $(1,1)\in
\C^2$ and so
\begin{equation}\label{grimprodbis}
(s_1,s_2)\in \R^2 \longrightarrow
(\gamma(s_1)+\gamma(s_2),\gamma(s_1)-\gamma(s_2))\in \c^2
\end{equation}
is a translating soliton for MCF with translating vector $(1,0)\in
\C^2$.

The translating solutions to mean curvature flow discovered by
Joyce, Lee and Tsui in \cite{JLT}, for the case $n=2$, are given
by $\mathcal{F}_t=\mathcal{F} + t(1,0)$ where $\mathcal F$ can be
described (see Section 1 in \cite{NT}) as follows: Let $w$ be a
curve in $\c$ whose curvature vector $\overrightarrow{k}$
satisfies $\overrightarrow{k}=w^\perp$. It can be chosen in such a
way that the angle $\theta $ that its tangent vector makes with
the $x$-axis has arbitrarily small oscillation. Then
\begin{equation}\label{exNT}
{\mathcal F}(x,y)=\left( \frac{|w(y)|^2-x^2}{2}-i\theta(y),x\,w(y) \right), \ (x,y)\in\r^2.
\end{equation}
It is still open the question posed in \cite{JLT} and \cite{NT}
about whether these translating solitons can arise as a blow-up of
a finite time singularity for Lagrangian mean curvature flow. It
would be very important to answer this question to develop a
regularity theory for the flow.


\section{New examples of Lagrangian translating solitons for MCF}

We start this section describing in the next Lemma a two-parameter
family of curves that provides a curious solution to the curve
shortening flow (CSF in short). Surprisingly some of them will be
the key ingredient for our construction of new examples of
Lagrangian translating solitons for MCF.

\begin{lemma}
Let $\alpha $ be a unit speed planar curve. Assume there exist
$a,b\in \r$, non null simultaneously, such that the curvature
function $\kappa_\alpha$ of $\alpha $ satisfies
\begin{equation}\label{kab}
\kappa_\alpha = a \langle \alpha, J \alpha' \rangle + b \langle
\alpha,  \alpha' \rangle
\end{equation}
where $'$ denotes derivative with respect to the arc parameter of
$\alpha$. Then the family of curves $\alpha_t=\sqrt{2at+1} \,
e^{i\frac{b}{2a}\log (2at+1)}\,\alpha$, with $2at+1>0$, is a
solution to the curve shortening flow
\begin{equation}\label{CSF}
\left( \frac{\partial}{\partial t} \alpha_t \right)^\perp =
\overrightarrow{\kappa_{\alpha_t}}
\end{equation}
such that $\alpha_0=\alpha$. Moreover, $\kappa_\alpha$ satisfies
the following o.d.e.
\begin{equation}\label{odeflux}
\kappa_\alpha \kappa_\alpha''-\kappa_\alpha'^2+
\kappa_\alpha^2(a+\kappa_\alpha^2)+b\,\kappa_\alpha'=0.
\end{equation}
\end{lemma}

\begin{remark}
{\rm In the limit cases $b=0$ and $a\rightarrow 0$ we recover well
known  solutions to the curve shortening flow:

If $b=0$, we have that the curvature vector of $\alpha $ verifies
$\overrightarrow{\kappa_\alpha} = a \, \alpha^\perp$ and so
$\alpha $ is a self-similar solution to CSF, contracting or
expanding according to $a<0$ or $a>0$ respectively; the flow
$\alpha_t=\sqrt{2at+1}\,\alpha$ is given by dilations of $\alpha$
in this case.

When $a \rightarrow 0$, we get now that
$\overrightarrow{\kappa_\alpha} = b (J\alpha)^\perp$
 and so $\alpha $ is a  spiral (see \cite{CZ}) solution to CSF with velocity $|b|$; the
flow $\alpha_t=e^{ibt}\,\alpha$ is given by rotations of $\alpha$
in this other case.
 }
\end{remark}

\begin{proof}
Using that the normal vector to $\alpha_t$ is given by
$J\alpha_t'/\sqrt{2at+1}$ and that
$\kappa_{\alpha_t}=\kappa_\alpha /\sqrt{2at+1}$, (\ref{CSF}) is
equivalent to $\langle \frac{\partial}{\partial t} \alpha_t, J
\alpha_t '  \rangle = \kappa_\alpha$. It is an exercise to check
that $\langle \frac{\partial}{\partial t} \alpha_t, J \alpha_t '
\rangle  = {\rm Im} \left( (a+ib) \overline{\alpha ' }\alpha
\right)$, which is precisely the condition satisfied by
$\kappa_\alpha$.

To prove the last part of the lemma, we define $f:=\langle \alpha
' , \alpha \rangle$ and $g:=\langle \alpha ' , J \alpha \rangle$
and so $\kappa_\alpha=b f- a g$. Using that $f'=1-\kappa_\alpha
g$, $g'=\kappa_\alpha f $ and $f^2+g^2=|\alpha|^2$, it is only a
long computation to check that $\kappa_\alpha$ satisfies
(\ref{odeflux}).
\end{proof}

In the next result, we make use of two families of planar curves
described in Lemma 1 (taking $a=\mp \cos \varphi$ and $b=\pm \sin
\varphi $ for a given $\varphi \in [0,\pi )$) in order to
construct many new Lagrangian translating solitons for MCF.

\begin{proposition}
Given $\varphi \in [0,\pi )$, let $\alpha=\alpha(t)$, $t\in I_1$,
and $\omega=\omega(s)$, $s\in I_2$, be unit speed planar curves
whose curvature vectors satisfy
\begin{equation}\label{cond1}
\overrightarrow{\kappa_\alpha}=-\cos\varphi\,\alpha^\perp+\sin\varphi\,(J\alpha)^\perp,
\
\overrightarrow{\kappa_\omega}=\cos\varphi\,\omega^\perp-\sin\varphi\,(J\omega)^\perp,
\end{equation}
where $\perp$ denotes normal component and $I_1$ and $I_2$ are
intervals of $\r$.

Let define $\alpha \ast \omega :I_1 \times I_2 \subset \r^2
\rightarrow \C^2$ by
\begin{eqnarray}\label{expl}
\varphi \neq \pi / 2: \ \  (\alpha \ast \omega)(t,s) =  \left(
\frac{|\omega(s)|^2 -|\alpha(t)|^2}{2\cos\varphi} \right.
\\ \left. +(\tan\varphi-i)(\arg\alpha'(t)+\arg\dot\omega(s)) \, ,
\, \alpha(t)\omega(s) \right)  \nonumber
\end{eqnarray}
and
\begin{eqnarray}\label{pi/2}
 \varphi = \pi / 2 : \ \ (\alpha \ast \omega)(t,s) =  \left(
\int_{t_0}^t \langle\alpha',J\alpha\rangle(x)\,dx \right.
\\ \left.
-\int_{s_0}^s \langle\dot\omega,J\omega\rangle(y)\,dy
-i(\arg\alpha'(t)+\arg\dot\omega(s)) \, , \,  \alpha(t)\omega(s)
\right) ,
 \nonumber
\end{eqnarray}
where ' and $\,\dot{}$ denote the derivatives respect to $t$ and
$s$ respectively, $t_0\in I_1$ and $s_0\in I_2$. Then  $\alpha
\ast \omega $ is a Lagrangian translating soliton for mean
curvature flow with translating vector $(1,0)\in \c^2$, whose
induced metric is $(|\alpha|^2 + |\omega|^2)(dt^2+ds^2)$ and its
Lagrangian angle map is $\arg\alpha'+\arg\dot\omega+\pi+\varphi$.
\end{proposition}

\begin{proof} The hypothesis on $\alpha $ and $\omega $ are
clearly equivalent to
$$\kappa_\alpha=\cos\varphi\,\langle\alpha',J\alpha\rangle+\sin\varphi\,\langle\alpha',\alpha\rangle, \
\kappa_\omega=-\cos\varphi\,\langle\dot\omega,J\omega\rangle-\sin\varphi\,\langle\dot\omega,\omega\rangle
$$
respectively. Then, looking at $\alpha$ and $\omega $ like complex
functions, (\ref{cond1}) is equivalent to
\begin{equation}\label{cond2}
\kappa_\alpha = {\rm Im}(e^{i\varphi}\alpha'\overline{\alpha}), \
\ \kappa_\omega = - {\rm
Im}(e^{i\varphi}\dot\omega\overline{\omega}).
\end{equation}
For any $t_0\in I_1$ and $s_0\in I_2$, using (\ref{cond1}) or
(\ref{cond2}), it is not difficult to check that the map $\alpha
\ast \omega $ can be written, up to a translation, in the
following common way for any $\varphi \in [0,\pi )$:
\begin{equation}\label{easy}
(\alpha \ast \omega)(t,s)=\left(e^{i\varphi}\left(\int_{s_0}^s
\dot\omega(y)\overline{\omega(y)}\,dy-\int_{t_0}^t
\alpha'(x)\overline{\alpha(x)}\,dx\right),
\alpha(t)\omega(s)\right)
\end{equation}
We denote $\Phi = \alpha \ast \omega$ and compute
$\Phi_t=\alpha'(-e^{i\varphi}\overline{\alpha},\omega)$ and
$\Phi_s=\dot\omega(e^{i\varphi}\overline{\omega},\alpha)$. Then we
obtain $|\Phi_t|^2=|\Phi_s|^2=|\alpha|^2 + |\omega|^2$ and
$(\Phi_t,\Phi_s)=0$. So $\alpha \ast \omega $ is a conformal
Lagrangian immersion whose induced metric is written as
$e^{2u}(dt^2+ds^2)$, with $e^{2u}=|\alpha|^2 + |\omega|^2$. So
$(t^*,s^*)$ is a singular point of $\alpha \ast \omega$ if and
only if $\alpha(t^*)=0=\omega(s^*)$.

Using that $e^{i\beta_{\alpha \ast \omega} }= e^{-2u} \det_\c
(\Phi_t,\Phi_s)$, it is not difficult to get that the Lagrangian
angle map $\beta_{\alpha \ast \omega}$ of $\alpha \ast \omega  $
is given by $\beta_{\alpha \ast \omega}
=\pi+\varphi+\arg\alpha'+\arg\dot\omega $. From (\ref{beta}) we
conclude that the mean curvature vector $H_{\alpha \ast \omega} $
of $\alpha \ast \omega $ is $ H_{\alpha \ast \omega}
=e^{-2u}\left(\kappa_\alpha J\Phi_t + \kappa_\omega J\Phi_s
\right) $.

On the other hand, $(1,0)^\perp = -e^{-2u} ({\rm
Im}(\Phi_t,(1,0))J\Phi_t+{\rm Im}(\Phi_s,(1,0))J\Phi_s)$ and hence
(\ref{cond2}) imply that $H_{\alpha \ast \omega} =(1,0)^\perp$.
\end{proof}

The conditions (\ref{cond1}) or (\ref{cond2}) are invariant by
rotations of the curves $\alpha$ and $\omega$. In the case
$\varphi =\pi/2$, $\alpha$ and $\omega$ must satisfy
$\overrightarrow{\kappa_\alpha}=(J\alpha)^\perp $ and $
\overrightarrow{\kappa_\omega}=-(J\omega)^\perp$. Thus, two
spirals $\alpha$ and $\omega$ with opposite velocities $\pm 1$
(see Remark 1)
 provide, under the construction $\alpha \ast \omega$, a Lagrangian
translating soliton for MCF. Since $\kappa_\alpha = \langle \alpha
, \alpha ' \rangle $ and $\kappa_\omega = - \langle \dot \omega,
\omega \rangle $, we get that the Lagrangian angle map in this
case $\varphi =\pi/2$ is given, up to a constant, by
$(|\alpha(t)|^2-|\omega(s)|^2)/2$.

In the same direction we now emphasize the case $\varphi =0$.

\begin{corollary}
Let $\alpha$ and $\omega$ self-similar solutions for the curve
shortening flow satisfying
$\overrightarrow{\kappa_\alpha}=-\alpha^\perp$ and
$\overrightarrow{\kappa_\omega}=\omega^\perp$. Then $\alpha \ast
\omega :I_1 \times I_2 \subset \r^2 \rightarrow \C^2$ given by
\begin{equation}\label{explss}
(\alpha \ast \omega)(t,s) =  \left( \frac{|\omega(s)|^2
-|\alpha(t)|^2}{2}  -i(\arg\alpha'(t)+\arg\dot\omega(s)),
\alpha(t)\omega(s) \right)
\end{equation}
is a Lagrangian translating soliton for mean curvature flow with
translating vector $(1,0)\in \c^2$.

By considering the straight lines $\alpha_0(t)=t$ and
$\omega_0(s)=s$, the circle $\alpha_1(t)=e^{it}$, joint to
self-shrinking curves $\alpha_{\mathcal S}$ and self-expanding
curves $\omega_{\mathcal E}$, we show up the following particular
examples:
\begin{itemize}
\item[(i)] $$(\alpha_0 \ast \omega_{\mathcal E})(t,s)
=\left(\frac{|\omega_{\mathcal
E}(s)|^2}{2}-i\arg\dot\omega_{\mathcal
E}(s)-\frac{t^2}{2},\,t\,\omega_{\mathcal E}(s)\right), $$ which
correspond to the Joyce, Lee and Tsui examples such as described
in (\ref{exNT});
\item[(ii)] $$(\alpha_1 \ast \omega_{\mathcal E})(t,s)
=\left(\frac{|\omega_{\mathcal
E}(s)|^2}{2}-i\arg\dot\omega_{\mathcal
E}(s)-it,\,e^{it}\omega_{\mathcal E}(s)\right),$$ for which
$\partial_t$ is a Killing vector field;
\item[(iii)] $$(\alpha_{\mathcal S} \ast \omega_0)(t,s)=
\left(\frac{s^2}{2}-\frac{|\alpha_{\mathcal
S}(t)|^2}{2}-i\arg\alpha_{\mathcal S}'(t),\alpha_{\mathcal
S}(t)s\right),$$ which satisfies that its Lagrangian angle map is
the angle that the tangent vector $\alpha_{\mathcal S}'(t)$ makes
with a fixed direction.
\end{itemize}
\end{corollary}

\begin{proof}
The result follows applying Proposition 2 with $\varphi =0$ and
taking into account that in the particular case (ii) the complete
induced metric is given by $(1+|\omega_{\mathcal
E}(s)|^2)(dt^2+ds^2)$ and in the particular case (iii) the
Lagrangian angle map is, up to a constant, the argument of
$\alpha_{\mathcal S}'(t)$.
\end{proof}

In Lemma 10.39 of \cite{CM} it is proved that any complete
self-shrinking planar curve is either a straight line passing
through the origin or it lies in a bounded set. The self-shrinking
curves found out by Abresh and Langer in \cite{AL} include a
countable family of non embedded closed curves. However, the
self-expanding planar curves $\omega_{\mathcal E}$ are embedded
and have two ends assymptotic to two straight lines (see for
example \cite{An} or \cite{EW}).

The totally geodesic Lagrangian plane is easily recovered in the
above construction by $(\alpha_0 \ast \omega_0)
(t,s)=\left(\frac{s^2-t^2}{2},t\,s\right)$. If we finally consider
the example $\alpha_1 \ast \omega_0$, we get the following result.
\begin{corollary}
Let define $\Phi:\r^2 \rightarrow \c^2$ by
$$ \Phi(t,s)=\left(\frac{s^2}{2}-it,e^{it}s\right).$$
Then $\Phi $ is a Hamiltonian stationary complete embedded
Lagrangian translating soliton for mean curvature flow with
translating vector $(1,0)\in \c^2$. In addition, $
\Phi(\r^2)=\mathcal{M} := \{ (z,w)\in \c^2 \, : \, w^2=2 {\rm Re}z
\, e^{-2i{\rm Im}z}, \,\, {\rm Re}z \geq 0 \}. $
\end{corollary}

\begin{proof}
We observe that $\Phi= \alpha_1 \ast \omega_0$. So, it is clear
that its induced metric is $(1+s^2)(dt^2+ds^2)$ and its Lagrangian
angle is $\beta(t)=3\pi/2+t$. Then $\Delta \beta =0$ and so $\Phi$
is Hamiltonian stationary.

Finally, it is clear that $\Phi(\r^2)\subset\mathcal{M}$. Given
$(z,w)\in \mathcal{M}$, we take $t=-{\rm Im}\,z$ and
$s=w\,e^{i{\rm Im}z}$. Since $s^2=2 {\rm Re}z\geq 0$, $s$ is well
defined and it is easy to check that $\Phi(t,s)=(z,w)$.
\end{proof}


\section{Classification of separable Lagrangian translating solitons}

In this section we characterize locally the examples of Lagrangian
translating solitons introduced in Proposition 2 under a
hypothesis that will allow us to separate variables in the
integration of the equations that translate (\ref{trl}).

\begin{theorem}\label{main}
Let $\phi: M^2 \rightarrow \c^2 $ be a Lagrangian translating
soliton for mean curvature flow with translating vector $\bf e$.
Assume that there exists a local isothermal coordinate $z=x+iy$
such that the smooth complex function $(\phi,{\bf e} )$ satisfies
$\frac{\partial^2}{\partial x \,
\partial y} (\phi,{\bf e} )=0$. Then $\phi $ is -up to
dilations- locally congruent to some of the following:
\begin{itemize}
\item[(i)] the product of a grim-reaper curve (\ref{grcurve}) and a straight line;
\item[(ii)] the product of two grim-reaper curves (see (\ref{grimprod}) or (\ref{grimprodbis}));
\item[(iii)] the example $\alpha \star \omega$ described in Proposition 2
for some $\varphi \in [0,\pi)$.
\end{itemize}
\end{theorem}

\begin{proof}
We start considering the translating vector ${\mathbf
e}=(1,0)\in\c^2$ without restriction and denoting $F=\langle \phi,
{\mathbf e} \rangle $ and $G =\langle \phi, J{\mathbf e} \rangle
$. Thus $\phi=(F+iG,\psi)$, where $\psi:M\rightarrow \C$ is the
second component of $\phi$. We will work in a local isothermal
coordinate $ z=x+iy$ on $M$ such that the induced metric, also
denoted by $\langle \,\, , \, \rangle $, is written as $\langle
,\rangle = e^{2u}|dz|^2$ with $|dz|^2$ the Euclidean metric. So we
have that
\begin{equation}\label{conf}
F_x^2+G_x^2+|\psi_x|^2=e^{2u}=F_y^2+G_y^2+|\psi_y|^2, \
F_xF_y+G_xG_y+\langle \psi_x,\psi_y \rangle =0
\end{equation}
and the Lagrangian character leads to
\begin{equation}\label{lagr}
F_yG_x-F_xG_y+\langle \psi_x,J\psi_y \rangle =0.
\end{equation}
Using (\ref{conf}) and (\ref{lagr}), taking into account that
$\psi_x$ and $\psi_y$ are both vectors in $\c$, it is not
difficult to get that
\begin{equation}\label{e2u}
e^{2u}=F_x^2+G_x^2+F_y^2+G_y^2, \ |\psi_x|^2=F_y^2+G_y^2, \
|\psi_y|^2=F_x^2+G_x^2 .
\end{equation}

From Proposition \ref{properties}, (\ref{beta}) and (\ref{trl}) we
also deduce that $F$ and $G$ must satisfy
\begin{equation}\label{trl1}
F_{xx}+F_{yy}=G_x^2+G_y^2, \quad G_{xx}+G_{yy}=-F_xG_x-F_yG_y
\end{equation}
and $\psi $ verifies
\begin{equation}\label{trl2}
\psi_{xx}+\psi_{yy}=-G_xJ\psi_x-G_yJ\psi_y.
\end{equation}

\vspace{0.3cm}

 From now on, by the hypothesis on separability, we can assume that the isothermal coordinate we are
 working in satisfies $(\phi,(1,0))_{xy}=0$. This means nothing
 but $F_{xy}=0=G_{xy}$. We remark that adding a constant to $F$ or $G$ produces a congruent immersion.
We make use of these two conditions in the
 following.

On the one hand, there exist smooth real functions $\xi=\xi(x)$
and $\theta=\theta(y)$ such that
$$ G(x,y)=-(\xi(x)+\theta(y)). $$
Then we consider planar curves $\alpha=\alpha(x)$, $x\in I_1
\subset \r$, and $\omega=\omega (y)$, $y\in I_2 \subset \r$, arc
length parameterized whose curvature functions are given by
$\kappa_\alpha(x)=\xi'(x)$ and $\kappa_\omega(y)=\dot\theta(y)$
respectively. Up to rotations, we can write
\begin{equation}\label{curves}
\alpha'(x)=e^{i\xi(x)}, \quad \dot\omega(y)=e^{i\theta(y)}
\end{equation}
and, up to a constant, we can also write
\begin{equation}\label{G}
G(x,y)=-\arg \alpha'(x)-\arg\dot\omega(y)=-\int\!\kappa_\alpha(x)dx -\int\!\kappa_\omega(y)dy.
\end{equation}
We also remark that, according to Proposition \ref{properties}(1),
the Lagrangian angle map $\beta $ of $\phi $ is given by
\[ \beta (x,y)= \xi(x)+\theta(y)+\beta_0, \ \beta_0\in \r . \]

 On the other hand, there exist smooth real
functions $A=A(x)$ and $B=B(y)$ such that
\begin{equation}\label{F}
 F(x,y)=A(x)+B(y).
\end{equation}
  Putting (\ref{F}) and (\ref{G}) in (\ref{trl1}), we can find
$\lambda, \, \mu \in\r$ such that $A$ and $B$  must satisfy the
following ordinary differential equations:
\begin{equation}\label{odesA}
\kappa_\alpha A' = \mu - \kappa_\alpha', \quad
A''=\kappa_\alpha^2-\lambda ,
\end{equation}
\begin{equation}\label{odesB}
\kappa_\omega \dot B = -\mu -\dot\kappa_\omega, \quad \ddot B =
\kappa_\omega^2+\lambda .
\end{equation}
We notice that the o.d.e.'s for $A$ and $B$ are the same
interchanging the pair $(\lambda,\mu)$ by $(-\lambda,-\mu)$. Let
us study (\ref{odesA}) for example. If $\kappa_\alpha \equiv 0$
then $\mu=0$ and $A(x)=-\lambda x^2/2 - b_1 x$, with $b_1 \in\r$,
up to a translation. If $\kappa_\alpha$ is non null, outside the
zeroes of $\kappa_\alpha $, we get $A(x)=-\log |\kappa_\alpha
(x)|+\mu\int dx / \kappa_\alpha(x)$, where $\kappa_\alpha$ is a
solution to
\begin{equation}\label{ode1}
\kappa_\alpha \kappa_\alpha''-\kappa_\alpha'^2+\mu
\kappa_\alpha'=\kappa_\alpha^2(\lambda-\kappa_\alpha^2).
\end{equation}
By the above observation, analogously if $\kappa_\omega \equiv 0$
then $\mu=0$ and $B(y)=\lambda y^2/2 +b_2 y$, with $b_2 \in\r$, up
to a translation. If $\kappa_\omega$ is non null, outside the
zeroes of $\kappa_\omega $, we get $B(y)=-\log |\kappa_\omega
(y)|-\mu\int dy / \kappa_\omega(y)$, where $\kappa_\omega$ is a
solution to
\begin{equation}\label{ode2}
\kappa_\omega\ddot\kappa_\omega-\dot\kappa_\omega^2-\mu
\dot\kappa_\omega=\kappa_\omega^2(-\lambda-\kappa_\omega^2).
\end{equation}
Hence we are devoted to study the o.d.e.'s (\ref{ode1}) and/or
(\ref{ode2}) in the following lemma, which deserves interest by
itself. We recognize (\ref{ode1}) and (\ref{ode2}) in Lemma 1
taking $a=\mp \lambda $ and $b=\pm \mu$ respectively.
\begin{lemma}
Given $\lambda,\mu\in \r$, consider the ordinary differential
equation
\begin{equation}\label{odek}
k\ddot k -\dot k^2 + k^2 (\lambda+k^2)=\mu \dot k .
\end{equation}
\begin{itemize}
\item If $(\lambda,\mu)=(0,0)$, then $\dot k^2 /k^2+ k^2 = \rho^2 \geq 0$ is a
first integral of (\ref{odek}) and $k(y)= \rho/\cosh(\rho y)$, $y\in\r$, is its solution satisfying
$\dot k (0)=0$.
\item If $(\lambda,\mu)\neq (0,0)$, let $k_w$ be the curvature of a unit speed planar curve
$w$ in $\c $ satisfying $k_w=-\lambda \langle \dot w, Jw \rangle -
\mu \langle \dot w, w \rangle$. Then $k_w$ is the general solution
of (\ref{odek}). Moreover, $k_w$ verifies:
\begin{enumerate}
\item $\frac{\textstyle (\dot k_w + \mu)^2}{\textstyle k_w^2}+k_w^2=(\lambda^2+\mu^2)|w|^2 $,
\item $ -\log |k_w| - \int\! \mu/k_w -i \int \! k_w =(\lambda + i \mu) \int \! \dot w \overline{w}
$,
\end{enumerate}
outside the zeroes of $k_w$.
\end{itemize}
\end{lemma}
{\it Proof of Lemma 2:\/} The case $(\lambda,\mu)=(0,0)$ is an
exercise. When $(\lambda,\mu)\neq (0,0)$, it was proved in Lemma 1
that $k_w$ satisfies (\ref{odek}). We define again $f:=\langle
\dot w, w \rangle$ and $g:=\langle \dot w, J w \rangle$ and so
$k_w=-\lambda g -\mu f$. Using that $\dot f=1-k_w g$, $\dot g=k_w
f $ and $f^2+g^2=|w|^2$, it is straightforward to check that $k_w$
satisfies part (1) in the Lemma. To prove part (2), we observe
that $ (-\log |k_w| - \int\! \mu/k_w -i \int \! k_w)\,\dot{}
=(\lambda + i \mu)(f+ig)= (\lambda + i \mu) \dot w \overline{w} $.
Finally, given arbitrary initial conditions $k_0=k(0)$ and
$k_1=\dot k(0)$ for (\ref{odek}), the system of equations
$-\lambda g(0)-\mu f(0)= k_0$, $\mu g(0)-\lambda f(0)=\mu +k_1$
has an unique solution since $(\lambda,\mu)\neq (0,0)$. This shows
that $k_w$ is the general solution of (\ref{odek}) and concludes
the proof of Lemma 2. \vspace{0.2cm}

We now proceed to integrate $\phi=(F+iG,\psi)$ collecting first the information from (\ref{F}),
(\ref{G}) and (\ref{e2u}). According to the above discussion, we must distinguish the following
cases:

{\bf Case (i):} $\kappa_\alpha \equiv 0 \equiv \kappa_\omega$. In particular $\mu=0$ and $G$ is
constant. Hence $\beta $ is constant too and so $\phi $ is minimal. Moreover, we have that
$$ F(x,y)=-\lambda x^2 /2 -b_1 x + \lambda y^2/2  +b_2
y $$ and $$ e^{2u(x,y)}=(\lambda x +b_1)^2+(\lambda y + b_2)^2 $$

{\bf Case (ii):} $\kappa_\alpha \equiv 0$ and $\kappa_\omega$ non
null. In particular $\mu=0$. We now get that
$$ F(x,y)=-\lambda x^2 /2 -b_1 x - \log |\kappa_\omega (y)|, \, G(y)=-\int\!\kappa_\omega(y)dy$$
and
$$ e^{2u(x,y)}=(\lambda x +b_1)^2+\dot \kappa_\omega(y)^2 /\kappa_\omega(y)^2+
\kappa_\omega(y)^2,
$$
where $\kappa_\omega$ is a solution of (\ref{ode2}) with $\mu=0$.

{\bf Case (iii):} $\kappa_\alpha $ non null and $\kappa_\omega \equiv 0$. In particular $\mu=0$.
Analogously we get that
$$ F(x,y)=- \log |\kappa_\alpha (x)|+\lambda y^2 /2 +b_2 y , \, G(x)=-\int \!\kappa_\alpha(x)dx$$
and
$$ e^{2u(x,y)}= \kappa_\alpha'(x)^2 /\kappa_\alpha(x)^2+ \kappa_\alpha(x)^2 +(\lambda y
+b_2)^2,
$$
where $\kappa_\alpha$ is a solution of (\ref{ode1}) with $\mu=0$.

{\bf Case (iv):} $\kappa_\alpha $ and $\kappa_\omega$ both non
null. We arrive at
$$ F(x,y)=- \log |\kappa_\alpha (x)|+\mu\int dx / \kappa_\alpha(x)-\log |\kappa_\omega (y)|-\mu\int dy
/ \kappa_\omega(y),$$ $$  G(x,y)=-\int\!\kappa_\alpha(x)dx -\int\!\kappa_\omega(y)dy$$ and
$$ e^{2u(x,y)}=(\kappa_\alpha'(x)-\mu)^2 /\kappa_\alpha(x)^2+
\kappa_\alpha(x)^2 + (\dot \kappa_\omega(y)+\mu)^2 /\kappa_\omega(y)^2+ \kappa_\omega(y)^2,$$ where
$\kappa_\alpha$ is a solution of (\ref{ode1}) and $\kappa_\omega$ is a solution of (\ref{ode2}).
 \vspace{0.3cm}

In order to use Lemma 2 we analyze the two given possibilities.
First we fix $\underline{(\lambda,\mu)\neq (0,0)}$. Using Lemma 2,
we know that $\omega $ and $\alpha $ must satisfy
$\kappa_\omega=-\lambda \langle \dot \omega, J\omega \rangle - \mu
\langle \dot \omega, \omega \rangle$ and $\kappa_\alpha=\lambda
\langle \alpha ' , J\alpha \rangle + \mu \langle \alpha ', \alpha
\rangle$ and, in addition, up to a constant we have that
\begin{equation}\label{FiG}
(F+iG)(x,y)=(\lambda + i \mu)\left(\int \dot \omega(y) \bar \omega(y) dy - \int \alpha'(x) \bar
\alpha (x) dx \right)
\end{equation}
and
\begin{equation}\label{e2ufin}
e^{2u(x,y)}=(\lambda^2 + \mu^2)(|\alpha(x)|^2+|\omega(y)|^2).
\end{equation}
In the cases (i), (ii) and (iii), necessarily $\lambda\neq 0$ since $\mu=0$. If we make changes of
parameters ($x\rightarrow x+b_1/\lambda$, $y\rightarrow y+b_2/\lambda$) then (\ref{FiG}) and
(\ref{e2ufin}) also hold (up to a translation) considering $\alpha(x)=x$ and $\omega(y)=y$ when
$\kappa_\alpha \equiv 0$ and $\kappa_\omega \equiv 0$ respectively.

Moreover, it is not difficult to get, taking into account (\ref{e2u}), (\ref{lagr}) and
(\ref{e2ufin}), that
\begin{equation}\label{psimod}
|\psi_x|^2=(\lambda^2 + \mu^2)|\omega|^2, \, |\psi_y|^2=(\lambda^2 + \mu^2)|\alpha|^2, \,
(\psi_x,\psi_y)=(\lambda^2 + \mu^2)\alpha' \bar \alpha \, \bar{\dot \omega}  \omega .
\end{equation}
Analyzing (\ref{trl2}) after considering (\ref{G}), using (\ref{curves}) and (\ref{psimod}), we
conclude that there exist two complex functions $c_i=c_i(x,y)$, $i=1,2$, such that
\begin{equation}\label{c1c2}
\psi_x=c_1 \alpha', \ \psi_y=c_2 \, \dot \omega, \quad (c_1)_x \,
\alpha'+(c_2)_y \, \dot \omega =0, \ \alpha c_1 = \omega c_2.
\end{equation}
Since $|c_1|^2=|\psi_x|^2$ and $|c_2|^2=|\psi_y|^2$, from (\ref{c1c2}) we can find two real
functions $\nu_i=\nu_i(x,y)$, $i=1,2$, in order to write
$c_1=\sqrt{\lambda^2+\mu^2}\,|\omega|e^{i\nu_1}$ and
$c_2=\sqrt{\lambda^2+\mu^2}\,|\alpha|e^{i\nu_2}$. The last two equations of (\ref{c1c2}) translate
into
\begin{equation}\label{nu1nu2}
 |\omega|   (\nu_1)_x \, \alpha'\,e^{i\nu_1} +|\alpha|  (\nu_2)_y \,\dot \omega \,  e^{i\nu_2}=0,
 \ |\omega| \alpha \,e^{i\nu_1} = |\alpha|\omega \, e^{i\nu_2},
\end{equation}
which lead to $(\nu_1)_x \, \omega \alpha' +(\nu_2)_y \, \dot \omega \alpha =0$. As $\alpha $ and
$\alpha '$ (resp.\ $\omega $ and $\dot \omega$) are necessarily linearly independent in this case,
we deduce that $(\nu_1)_x=0=(\nu_2)_y$ and hence there is a constant $\nu_0$ such that
$|\omega|e^{i\nu_1}/\omega=|\alpha|e^{i\nu_2}/\alpha=e^{i\nu_0}$ thanks to the last equation in
(\ref{nu1nu2}). Using the first two equations of (\ref{c1c2}), we arrive at
$\psi_x=\sqrt{\lambda^2+\mu^2}\,e^{i\nu_0}\alpha'\omega$ and $\psi_y=\sqrt{\lambda^2+\mu^2}\,
e^{i\nu_0}\alpha\dot\omega$. Thus, up to a rotation and a translation, we finally get that
\begin{equation}\label{psifin}
\psi(x,y)=\sqrt{\lambda^2+\mu^2}\,\alpha(x)\omega(y).
\end{equation}
Therefore we conclude from (\ref{FiG}) and (\ref{psifin}) that
$$
\phi(x,y)=\left((\lambda +i\mu)\left(\int \!\dot \omega(y) \bar \omega(y) dy - \int\! \alpha'(x)
\bar \alpha (x) dx \right),\sqrt{\lambda^2+\mu^2}\,\alpha(x)\omega(y) \right)
$$
where $\alpha $ and $\omega$ satisfy $\kappa_\omega=-\lambda
\langle \dot \omega, J\omega \rangle - \mu \langle \dot \omega,
\omega \rangle$ and $\kappa_\alpha=\lambda \langle \alpha ' ,
J\alpha \rangle + \mu \langle \alpha ', \alpha \rangle$. Up to
dilations, there is no restriction taking $\lambda +i\mu
=e^{i\varphi}$, with $\varphi\in [0,2\pi)$. So this is exactly the
common expression (\ref{easy}) for the examples $\alpha \star
\omega $ introduced in Proposition 2. Interchanging the roles of
$\alpha$ and $\omega$, it is enough to consider $\varphi\in
[0,\pi)$. The conclusion is that $\phi$ is one of the examples
mentioned in part (iii) of the statement of Theorem 1.
\vspace{0.2cm}

 We finally study the remaining case $\underline{(\lambda,\mu)= (0,0)}$. We remark that in cases (i), (ii) and
 (iii) we only have to consider $\lambda=0$ since $\mu$ was necessary zero there.

 In case (i), $\lambda = 0$ implies that $u$ is constant and so the immersion is flat besides
 minimal. Thus it is totally geodesic. Recall that $\alpha_0 \star
\omega_0 $ recovers a totally geodesic Lagrangian plane.

 In case (ii), following Lemma 2, up to a constant we get that
 $$
(F+iG)(x,y)=-b_1 x + \log \cosh (\rho y) - i \rho \int \frac{dy}{\cosh \rho y}.
 $$
In the coordinates $(t,s)=-\rho(x,y)$ and putting $-b_1/\rho=\sinh \delta $, $\delta \in \R$, we
rewrite
$$
(F+iG)(t,s)=-\sinh \delta \,  t + \log \cosh s + i  \int \frac{ds}{\cosh s} = -\sinh \delta \, t +
\gamma (s),
 $$
where $\gamma (s) $ (see (\ref{grcurve})) is just the graph
$(-\log v,v)$, $v\in (-\pi/2,\pi/2)$, parameterized by arc length.
A similar study of (\ref{trl2}) like in the previous case using
now in (\ref{conf}), (\ref{lagr}) and (\ref{e2u}) the above
expressions of $F$ and $G$ leads to
$$
\psi(t,s)=t+\sinh \delta \, \gamma (s).
$$
Then we get that $\phi (t,s)= A (\gamma (s), t)$, where $A$ is the
matrix $\left( \begin{array}{cc}  1 & -\sinh \delta \\
\sinh \delta & 1 \end{array} \right)$. Thus we arrive at (i) in
the statement of Theorem 1.

Case (iii) is completely analogous to case (ii) changing $\omega $
by $\alpha $ and $b_1$ by $-b_2$ so that we get the same
conclusion.

In case (iv), applying twice Lemma 2 and the same argument that in
case (ii), we deduce that
 $$
(F+iG)(s_1,s_2)= \gamma (s_1) + \gamma (s_2), \quad
\psi(s_1,s_2)=\gamma (s_1) - \gamma (s_2).
 $$
Hence we arrive at (ii) in the statement of Theorem 1.
\end{proof}

\begin{corollary}
Let $\phi :M \rightarrow \C^2 $ be a Hamiltonian stationary (non
totally geodesic) Lagrangian translating soliton for mean
curvature flow. Then $\phi(M) $ is -up to dilations- an open
subset of the Lagrangian $\mathcal{M}$ given in Corollary 2.
\end{corollary}

\begin{proof} We follow the same election for the translating
vector and use the same notation that at the beginning of the
proof of Theorem 1.  We can associate to any Lagrangian immersion
$\phi :M \longrightarrow \c ^2$ a differential form $\Upsilon $ on
$M$ (see \cite{CU1}) defined by
\[ \Upsilon (z)=\bar{h}(z)dz, {\rm \ with \ } h(z)= \omega (
\partial_{\bar z}, H ), \]
where $ z=x+iy$ is a local isothermal coordinate on $M$ and
$\omega $ is extended $\c$-linearly to the complexified tangent
bundles. Then (\ref{beta}) translates into $h=\beta_{\bar z}$,
with $\beta $ the Lagrangian angle map of $\phi$, and the Coddazi
equation of $\phi$ gives (see \cite{CU1}) ${\rm Im}(h_z)  =  0 $.

Thus $\bar{h}_{\bar{z}}=h_z=\beta_{z \overline{z}}=0$ since
$\beta$ is harmonic because $\phi$ is Hamiltonian stationary.
Hence $\Upsilon$ is a holomorphic differential and we can
normalize $h\equiv -1/2$.

Using Proposition 1.(1), we have that also $h=-G_{\bar z}$ and so
$G_x\equiv -1$ and $G_y\equiv 0$ after the above normalization. In
particular, $G_{xy}=0$. Looking at the second equation of
(\ref{trl1}) we easily deduce that $F_x=0$ and then $F_{xy}=0$. We
have proved that $\phi $ verifies the hypothesis of Theorem 1 and
necessarily must be one the examples $\alpha \star \omega$
associated to certain $\varphi \in [0,\pi) $. We know from
Proposition 2 that its induced metric is conformal and, using the
expression of its Lagrangian angle map, we get that $\alpha \star
\omega$ is Hamiltonian stationary if and only if $\kappa_\alpha '
+ \dot \kappa_\omega = 0$. Using (\ref{odeflux}) we obtain that
$\kappa_\alpha \equiv c_1 \in\r$ and $ \kappa_\omega \equiv c_2
\in \r$ such that
\begin{equation}\label{hsl}
c_1^2(c_1^2-\cos\varphi)=0=c_2^2(c_2^2+\cos \varphi).
\end{equation}
If $c_1=0$, $\alpha$ must be a line and this implies that $\varphi
=0 $ and, following the notation of Corollary 1,
$\alpha=\alpha_0$. Using now (\ref{hsl}) we have that $c_2=0$ and
a similar reasoning gives that $\omega = \omega_0$. In this case,
$\phi $ corresponds to a totally geodesic Lagrangian plane.

And if $c_1\neq 0$, from (\ref{hsl}) it follows that
$c_1^2=\cos\varphi$, $0\leq \varphi < \pi/2$, and $c_2=0$. This
last implies that $\omega$ must be a line, $\varphi =0$ and
$\omega = \omega_0$. Thus $c_1=1$ and we finally deduce that
$\alpha =\alpha_1$. Therefore we arrive at the example
$\Phi=\alpha_1 \star \omega_0$ and Corollary 2 finishes the proof.
\end{proof}

\end{document}